\documentclass[11pt,twoside]{article}
\usepackage{latexsym,epsfig,graphpap,graphics,amsmath,amssymb,float}
\usepackage[utf8x]{inputenc}
\newtheorem{thm}{Theorem}[section]

\newtheorem{prop}[thm]{Proposition}
\newtheorem{defn}[thm]{Definition}
\newtheorem{remark}[thm]{Remark}

\newtheorem{example}[thm]{Example}
\numberwithin{equation}{section}

\usepackage{plain}
\begin{document}
\title{Charnes-Cooper Scalarization and Convex Vector Optimization}
\author{ Poonam Kesarwani\footnote{Department of Mathematics and Statistics, Indian Institute of Technology Kanpur, India} and Joydeep Dutta\footnote{Department of Economics Sciences, Indian Institute of Technology Kanpur, India} }
\maketitle

\begin{abstract}
Our aim in this article is two-fold. We use the Charnes-Cooper scalarization technique to develop KKT type conditions to completely characterize Pareto minimizers of convex vector optimization problems and further, we use that scalarization technique to develop a simple and efficient algorithm for convex vector optimization problems. Numerical examples are presented to illustrate the use of our algorithm.
\end{abstract}

\section{Introduction}
In this paper, we shall focus on the well known convex vector optimization problem given below, which we label as CVOP:
\begin{eqnarray*}
&&\min f(x):=(f_1(x),\ldots,f_m(x)),\\
&&{\rm subject~to}~~ x\in X,
\end{eqnarray*}
where each $f_i:\mathbb{R}^n\rightarrow \mathbb{R}$ is a convex function and $X$ is a closed convex subset of $\mathbb{R}^n$. Let us denote the set $I:=\{1,2,\ldots, m\}$ for the sake of convenience.  The ``min'' in the above problem is taken in the vector sense where a partial order is induced in the image space $\mathbb{R}^m$, by the non-negative cone $\mathbb{R}^m_+$. The partial ordering says that $x\ge y$, if $x-y\in \mathbb{R}^m_+$, which can equivalently be written as $x_i\ge y_i$ for $i\in I$, where $x_i$ and $y_i$ represents the $i$-th component of the vectors $x$ and $y$ respectively.\\
Two Types of solution concepts are popular in the literature of vector optimization problems. These are namely a Pareto minimizer and a weak Pareto minimizer. A vector $x^*\in X$ is called a Pareto minimizer of CVOP if there exists no ${x}\in X$ such that $f_i({x})\le f_i(x^*)$ for all $i\in I$, and, $f_r({x})\le f_r(x^*)$ for some $r\in I$.
More compactly, it means that there exists no ${x}\in X$ such that 
$$f(x)-f(x^*)\in -\mathbb{R}^m_+\setminus \{0\}.$$
On the other hand, a vector $x^*\in X$ is called a weak Pareto minimizer if there exists no ${x}\in X$ such that $f_i({x})< f_i(x^*)$ for all $i\in I$, which equivalently be written as $f(x)-f(x^*)\in -\text{int}(\mathbb{R}^m_+)\setminus \{0\}.$ Though from a practical point of view, the notion of a Pareto minimizer is more useful than its weak counterpart. However, the notion of the weak Pareto minimizer is mathematically more tractable, and thus, one will find a huge literature on these solutions. Note that every Pareto minimizer is a weak Pareto minimizer of the problem and the converse need not be true. The image of the set of Pareto minimizers of CVOP is called the efficient frontier, and it has been shown, for example, in Ehrgott \cite{EHR} that the efficient frontier lies on the boundary of the image $f(X)$ of the feasible set $X$ under the objective function $f$. \\
One of the key approach to study vector optimization problem is by relating a scalar optimization problem to the given vector optimization problem. This method is called scalarization and plays a key role in developing a coherent theory and effective numerical algorithms for the vector optimization problem (see Chankong and Haimes \cite{chankong2008multiobjective} for more details).\\
As far as CVOP is concerned, the following result completely characterizes the weak minimizers through a simple scalarization technique called the weighted sum scalarization.
\begin{prop}\label{prop1}
Consider the problem CVOP and consider weights $\tau_1\ge 0$,$\tau_2\ge 0$,$\ldots, \tau_m\ge 0$, and $x^*\in X$ such that $x^*$ solves the problem (CP$_\tau$),
\begin{equation}\label{eq11}
\min\limits_{x\in X}(\tau_1f_1(x)+\ldots+\tau_mf_m(x)),
\end{equation}
then $x^*$ is a weak Pareto minimizer of CVOP. Conversely, if $x^*\in X$ is a weak Pareto minimizer of CVOP, then there exist scalars $\tau_1\ge 0$,$\ldots, \tau_m\ge 0$, all of which are not simultaneously zero, such that $x^*$ is a minimizer of CP$_\tau$.
\end{prop}
This result immediately provides us a KKT type necessary and sufficient condition which completely characterizes a weak minimizer of CVOP. In fact, by using Proposition \ref{prop1}, we can conclude that 
a vector $x^*\in X$ is a weak minimizer of CVOP if and only if there exists $0\not= \tau \in \mathbb{R}^m_+$ such that 
\begin{equation}\label{eq12}
0\in \sum\limits_{i\in I}\tau_i \partial f_i(x^*)+N_X(x^*),
\end{equation}
where $\partial f_i(x^*)$ is the subdifferential at $x^*$ of $f$ and $N_X(x^*)$ denotes the normal cone of $X$ at $x^*$. For the definition of the subdifferential of a convex function and the normal cone to a convex set at a given point see Rockafellar \cite{rockafellar2015convex}. Further, if each $f_i$ is differentiable, then we can replace $\partial f_i(x^*)$ by $\nabla f_i(x^*)$ in \eqref{eq12}. Moreover, if $X$ is given by convex inequality constraints, \textit{i.e.},
\begin{equation}\label{eq13}
X=\{x\in \mathbb{R}^n:g_j(x)\le 0,~\text{for all }j=1,\ldots,k\},
\end{equation}
where each $g_j:\mathbb{R}^n\rightarrow \mathbb{R}$ is a convex function, for all $j=1,\ldots,k$, and, the Slater condition holds, \textit{i.e.}, there exists $\hat{x}\in \mathbb{R}^n$ such that $g_j(\hat{x})<0$ for all $j=1,\ldots,k$, then we have (see Rockafellar \cite{rockafellar2015convex}),
\begin{equation*}
N_X(x^*)=\{v\in \mathbb{R}^n: v=\sum\limits_{j\in K} \lambda_j \xi^j,\lambda_j\ge 0, \xi^j\in \partial g_j(x^*),\lambda_jg_j(x^*)=0, \forall~j\in K\},
\end{equation*}
where the set $K:=\{1,2,\ldots,k\}$. This immediately allows us to state a KKT type optimality condition for CVOP, which we do through the following theorem.

\begin{thm}\label{t1}
Consider the problem CVOP, with $X$ described by \eqref{eq13}. Let the Slater condition holds. Then $x^*\in X$ is a weak Pareto minimizer of CVOP if and only if there exists vectors $\tau \in \mathbb{R}^m_+$, $\tau\ne 0$ and $\lambda\in \mathbb{R}^k_+$ such that
\begin{description}

\item $(i)$ $0\in \sum\limits_{i\in I} \tau_i\partial f_i(x^*)+\sum\limits_{j\in K} \lambda_j\partial g_j(x^*)$,
\item $(ii)$ $\lambda_jg_j(x^*)=0,$ for all $j\in K$.

\end{description}

\end{thm}
For a more detailed analysis of the above results and their proofs, the reader is referred to the monographs by Jahn \cite{JAHN}, Ehrgott \cite{EHR}, Chankong \& Haimes \cite{chankong2008multiobjective}, Miettinen \cite{MIET} and the references therein.\\
Thus, Theorem \ref{t1} shows that it is possible to completely characterize a weak Pareto minimizer of CVOP through KKT conditions. Now we ask ourself a question: Can this type of result is true for a Pareto minimizer, \textit{i.e.}, can we completely characterize a Pareto minimizer using KKT type optimality conditions? \\
Developing a necessary condition for a Pareto minimizer is indeed simple. By observing that a Pareto minimizer is a weak Pareto minimizer, Theorem \ref{t1} guarantees the existence of a vector $\tau \in \mathbb{R}^m_+$, $\tau\ne 0$ and $\lambda\in \mathbb{R}^k_+$ such that $(i)$ and $(ii)$ hold. The problem comes when we are considering the sufficiency of the conditions. Unfortunately, the conditions $(i)$ and $(ii)$ can only show that $x^*$ is a weak Pareto minimizer of CVOP and no way, it can guarantee that $x^*$ is a Pareto minimizer. However, in $(i)$ if we have $\tau \in \text{int}(\mathbb{R}^m_+)$, then we can show that $x^*$ is a Pareto minimizer. Therefore, in order to completely characterize a Pareto minimizer for CVOP through optimality condition, we need to develop a necessary optimality condition in which we must have $\tau \in \text{int}(\mathbb{R}^m_+)$. For this, one either needs the Pareto minimizer to have some additional properties or some regularity conditions are needed involving the objective functions and the constraint functions. When $\tau \in \text{int}(\mathbb{R}^m_+)$ in the condition $(i)$, the KKT conditions are called the \textbf{strong KKT} conditions. Historically, it was Kuhn and Tucker in their seminal paper \cite{kuhn1951proceedings} of 1951 develops a notion of KT- proper Pareto solutions for which the strong KKT holds. To prove the strong KKT conditions Kuhn-Tucker used a regularity condition called Kuhn-Tucker 
constraint qualification \cite{kuhn1951proceedings}. Later in 1994, Maeda \cite{maeda1994constraint} showed that under Guignard type regularity condition a Pareto minimizer becomes a KT-proper Pareto minimizer. However, the simplest approach to strong KKT conditions can be shown using the notion of Geoffrion proper minimizers, a notion developed by Geoffrion \cite{GEO} in 1968.
\begin{defn}(Geoffrion \cite{GEO})
A vector $x^*\in X$ is called a Geoffrion proper minimizer of CVOP if $x^*$ is a Pareto minimizer and if there exists a number $M>0$ such that for all $i\in I$ and $x\in X$ satisfying $f_i(x)< f_i(x^*)$, there exists a $j\in  I$ such that $f_j(x^*)< f_j(x)$ and
\begin{eqnarray}
\frac{f_i(x^*)- f_i(x)}{f_j(x)- f_j(x^*)}\leq M.
\end{eqnarray}
\end{defn}
Geoffrion \cite{GEO} shows that if $x^*$ is Geoffrion proper minimizer of CVOP with smooth input data and the Kuhn-Tucker 
constraint qualification \cite{kuhn1951proceedings} is satisfied at $x^*$, then there exists vectors $\tau \in \text{int}(\mathbb{R}^m_+)$, and $\lambda\in \mathbb{R}^k_+$ such that
\begin{itemize}
\item  $ \sum\limits_{i\in I} \tau_i\nabla f_i(x^*)+\sum\limits_{j\in K} \lambda_j\nabla g_j(x^*)=0$,
\item  $\lambda_jg_j(x^*)=0,$ for all $j\in K$.
\end{itemize}
If we replace the Kuhn-Tucker Constraint qualification by the Slater condition in the above result, we will get the same strong KKT condition for Geoffrion proper minimizer. We would like to mention here that the strong KKT conditions for a Geoffrion proper minimizer also holds even if some of the objective functions are not differentiable. For such functions, we just have to replace the gradient with the subdifferential of the functions. Therefore, Geoffrion proper minimizer gives rise to strong KKT type optimality conditions.  It is well known that every Pareto minimizer is not a Geoffrion proper minimizer (see Example 2.48 of \cite{EHR}). Hence, we are now asking whether it is possible to develop a strong KKT type optimality conditions for Pareto minimizers. \\
In Section 2, we shall show that this could be achieved by a combination of Charnes and Cooper scalarization \cite{charnes1961management} and Abadie type regularity condition which we shall describe at the beginning of the section. In Section 3, we shall show that Charnes and Cooper scalarization can be efficiently used to develop an algorithm to solve convex vector optimization problems where we need the feasible set to be closed and $f_i$'s to be a proper lower semi-continuous function for all $i\in I$. We shall also illustrate our algorithm through numerical examples. 
\section{Strong KKT conditions for Pareto minimizers}
There two main purposes for which the Charnes-Cooper scalarization is introduced. The first is to completely characterize Pareto minimizer through KKT conditions. This is what we do in this section by showing how the Charnes-Cooper scalarization scheme can be used to derive strong KKT condition. The second purpose is to develop an algorithm to generate Pareto minimizer of CVOP which is done in the next section.
 
 Consider an arbitrary feasible point $x_0\in X$ of CVOP and define the following scalar problem $P(f,X,x_0) $ (\cite{ehrgott2005saddle}): 
 \begin{equation} 
 \left\{ \begin{array}{ll} \min \sum\limits_{i\in I}f_i(x), & \label{c1e1}\\
{ \text{subject to }}~~ f_i(x)-f_i(x_0)\leq 0, &\text{ for all }  i\in I, \\  
x\in X. & \end{array}\right.
\end{equation}
The advantage of the above scalar problem is that it provides a one-one correspondence between Pareto minimizer of CVOP and optimal solutions of the problem $P(f,X,x_0)$.  
\begin{thm}[Ehrgott\cite{ehrgott2005saddle}]\label{c1tm34}
A point $x_0\in X$ is a Pareto minimizer of CVOP if and only if $x_0$ is an optimal solution of the problem $P(f,X,x_0)$.  
\end{thm}
Observe that even if we consider the objective functions to be non-convex, Theorem \ref{c1tm34} is true, in other words, the characterization of a Pareto-minimizer through Charnes-Cooper scalarization still holds. However, computing the global minimizer of a non-convex optimization problem is hard, and we can only assure the computation of local Pareto minimizers which are not of special interest. So, we shall focus only on convex problems.\\
 Note that if $x_0$ is a Pareto minimizer of CVOP, then Slater constraint qualification cannot be satisfied for the problem $P(f,X,x_0)$.  Therefore, if we want to use the scalar problem $P(f,X,x_0)$ to figure out the necessary and optimality condition for Pareto minimizer, we need some regularity condition to be satisfied. In this direction, we define the following Abadie type regularity condition for CVOP which helps us in further investigation.

\begin{defn}
Consider CVOP with feasible set $X$ and for given $x_0\in X$, define the set $\hat{X}(x_0):=\{x\in X:f_i(x)\leq f_i(x_0),~i\in I\}.$ Then, we say that the problem satisfies the strong Abadie type regularity condition at a point $x_0\in X$ if the tangent cone $T_{\hat{X}}(x_0)$ of the set $\hat{X}$ at the point $x_0$ equals
\begin{equation}\label{eqN122}
V(x_0):=\{h\in \mathbb{R}^n: f_i^{'}(x_0,h)\leq 0,~ g_k^{'}(x_0,h) \leq 0, \text{ for all } i\in I\text{ and }k\in \mathcal{R}(x_0)\}.
\end{equation}
\end{defn}
\noindent Recall that $\mathcal{R}(x_0)$ is the active constraint set, \textit{i.e.}, $\mathcal{R}(x_0)=\{k\in K: g_k(x_0)=0\}$ and and the tangent cone $T_X(x_0):=\{d\in \mathbb{R}^n: \exists ~\{t_k\}\subseteq    \mathbb{R}_+, \{y_k\}\subseteq   X~\text{ subject to }~ y_k\rightarrow x_0,~ t_k(y_k-x_0)\rightarrow d\}.$ \\
Note that the feasible set of  $(f,X,x)$ is ${\hat{X}}(x_0)$. One can observe that the strong Abadie type regularity condition is same as the Abadie constraint qualification (see \cite{dhara2011optimality}) for the scalar problem $P(f,X,x_0)$.  Before we derive the KKT conditions for Pareto minimizer, let us establish that the problem CVOP satisfies the strong Abadie type regularity condition through two examples in which one is for smooth CVOP and other is for non-smooth CVOP.

\begin{example}(Chandra et al.\cite{chandra2004regularity}) Consider the bi-objective optimization problem, where $f=(f_1,f_2)$ with $f_1(x)=x^2,~f_2(x)=x~{ and}~X=\{x\in \mathbb{R}:g(x)\leq 0\}$, where $g(x)=-x$. It is clear that $x_0=0$ is a Pareto minimizer of the problem. Here, the set $\hat{X}(x_0)=\{x\in \mathbb{R}: x^2\leq 0, x\leq 0, -x\leq 0\}=\{0\}$ and it is clear that $T_{\hat{X}(x_0)}(0)=\{0\}$.  Further, since the constraint function $g$ is active at $x_0=0$ and $f_1^{'}(x_0,h)=0 $, the set
\begin{equation*}
V(x_0) =\{h\in \mathbb{R}: h\leq 0,-h\leq 0 \}=\{0\}.
\end{equation*}
 Hence, $T_{\hat{X}(x_0)}=V(x_0)$, which implies that considered problem satisfies the strong Abadie type regularity condition at the point $x_0=0$.
 \end{example}
 \begin{example}
 Consider the bi-objective optimization problem  with $f=(f_1,f_2)$ where $f_1(x)=x,~f_2(x)=|x|~{ and}~X=\mathbb{R}.$ It is clear that $x_0=0$ is a Pareto minimizer of the problem. Here, the set $\hat{X}(x_0)=\{x\in \mathbb{R}: x\leq 0, |x|\leq 0\}=\{0\}$ and hence, $T_{\hat{X}(x_0)}(0)=\{0\}$. Now form \eqref{eqN122}, the set
$V(x_0) =\{h\in \mathbb{R}:h\leq 0, |h|\leq 0 \}=\{0\},$ which implies that $T_{\hat{X}(x_0)}=V(x_0)$. Hence,the strong Abadie type regularity condition holds at the point $x_0=0$.  
\end{example}
The idea of using the term strong Abadie type regularity is motivated by the fact that all the objective functions of CVOP appear in the formulation of the Abadie type regularity condition. For a weaker version of the Abadie type regularity condition for a multiobjective optimization problem see \cite{chandra2004regularity}.
\begin{thm}
Consider CVOP in which all objective functions $f_i$'s and constraint functions $g_r$'s are smooth. Let $x_0$ be a Pareto minimizer of CVOP and the strong Abadie type regularity condition be satisfied at $x_0$. Then there exists $\lambda^*\in \mathbb{R}^m_+$ and $\mu^*\in \mathbb{R}^k_+$ such that
\begin{eqnarray}
 \sum\limits_{i\in I}\lambda_i^*\nabla f_i(x_0)+\sum\limits_{j\in K}\mu_j^*\nabla g_j(x_0)=0, \label{eqN1}\\
 \mu_j^*g_j(x_0)=0, \text{ for all } j\in K,\label{eqN2}\\
 \lambda_i^*>0, \text{ for all }i\in I\label{eqN3}.
\end{eqnarray}
Conversely, if there exists $\lambda^*\in \mathbb{R}^m_+$ and $\mu^*\in \mathbb{R}^k_+$ satisfying  \eqref{eqN1}-\eqref{eqN3}. Then $x_0$ is a Pareto minimizer of CVOP. 
\end{thm}

\textit{Proof:} Let $x_0$ be a Pareto minimizer of CVOP. Then, using Theorem~\ref{c1tm34},  $x_0$ is an optimal solution of the problem $P(f,X,x_0)$. Hence, from the standard necessary optimiality conditions, we know that
\begin{equation}\label{eqN125}
\left\langle \sum\limits_{i\in I}\nabla f_i(x_0),d\right\rangle\geq 0, \text{ for all }d\in T_{\hat{X}(x_0)}(x_0). 
\end{equation}
Since, $\hat{X}(x_0)$ is a convex set, the tangent cone $T_{\hat{X}(x_0)}(x_0)$ is a closed convex cone. Using the fact that the strong Abadie type regularity condition satisfies at $x_0$, \eqref{eqN125} becomes,
\begin{eqnarray}\label{eqN126}
-\sum\limits_{i\in I}\nabla f_i(x_0)\in (T_{\hat{X}(x_0)}(x_0))^\circ=(V(x_0))^\circ
\end{eqnarray}
Using \eqref{eqN122} and the definition of polar cone, we have
\begin{equation}\label{eqN127}
(V(x_0))^\circ=\{\sum\limits_{i\in I}\lambda_i\nabla f_i(x_0)+\sum\limits_{r\in K}\mu_r\nabla g_r(x_0): \lambda\in \mathbb{R}^m_+, \mu\in \mathbb{R}^k_+, \mu_rg_r(x_0)=0, \text{ for all }r\in K\}.
\end{equation}  Therefore, \eqref{eqN126} says that there exists $\lambda\in \mathbb{R}^m_+, \mu\in \mathbb{R}^k_+$ such that $\mu_rg_r(x_0)=0$ for all $r\in K$ and $$0=\sum\limits_{i\in I}\nabla f_i(x_0)+\sum\limits_{i\in I}\lambda_i\nabla f_i(x_0)+\sum\limits_{r\in K}\mu_r\nabla g_r(x_0).$$
Further, the above conditions can be rewritten as 
\begin{enumerate}
\item $0=\sum\limits_{i\in I}\lambda^*_i\nabla f_i(x_0)+\sum\limits_{r\in K}\mu^*_r\nabla g_r(x_0),$
\item $\mu_r^*g_r(x_0)$ for all $r\in K$,
\item $\lambda^*_i>0$, for all $i\in I$.
\end{enumerate}
where $\lambda^*_i=\lambda_i+1$ for all $i\in I$ and $\mu^*_r=\mu_r$ for all $r\in K$, which completes the necessary part.\\
For sufficient part, let there exists $\lambda^*\in \mathbb{R}^m_+$ and $\mu^*\in \mathbb{R}^k_+$ satisfying  \eqref{eqN1}-\eqref{eqN3}. Then from the Theorem 3.26 of \cite{EHR}, $x_0$ is a Geoffrion proper minimizer of CVOP. Hence $x_0$ is a Pareto minimizer of CVOP which completes the proof.\hfill$\Box$

\begin{remark} 
When input data is non-smooth, then we only get necessary optimality condition which is not same \eqref{eqN1}-\eqref{eqN3}. Since, for non-smooth case, \eqref{eqN127} turns out to be
\begin{multline}\label{eqN129}
(V(x_0))^\circ=cl\{\sum\limits_{i\in I}\lambda_iu_i+\sum\limits_{r\in K}\mu_rv_r :  \lambda\in \mathbb{R}^m_+,\mu\in \mathbb{R}^k_+,u_i\in \partial f_i(x_0),\\ v_r\in \partial g_r(x_0), \mu_rg_r(x_0)=0, \text{ for all }i\in I\text{ and }r\in K\},
\end{multline}
Hence, the necessary condition for $x_0$ to be a Pareto minimizer under the strong Abadie type regularity condition is given by
$$0\in \sum\limits_{i\in I}\partial f_i(x_0)+(V(x_0))^\circ,$$
where $V(x_0)$ is given by \eqref{eqN129}. The verification of the strong Abadie type regularity condition is not easy which makes Pareto minimizer less algorithm friendly and, almost all algorithms focus on computing weak Pareto minimizer. 

\end{remark}

\section{An Algorithm for computing Pareto minimizers}
It is well known that there are lots of scalarization technique to solve CVOP consider for example, weighted sum scalarization, $\epsilon$-constraint method, Benson method and many more. The way by which Charnes and Cooper scalarization is different form others is that there is no need of any parameters (\textit{e.g.} weighted sum scalarization needs weights and  $\epsilon \in \mathbb{R}_+^m$-constraint method needs $\epsilon$ etc.) when we use this scalarization to compute Pareto minimizers in the algorithm. Hence, using Charnes and Cooper scalarization for computing Pareto minimizers is more suitable rather than using other scalarization techniques. Let us discuss the algorithm to find Pareto minimizers using Charnes and Cooper scalarization and we shall call the proposed algorithm as CC1. In the algorithm below, the set $Sol(P_k)$, denotes the solution set of the scalar program $P_k$.
\begin{description}
\item \textit{Step 1.} (Initialization) Choose $x_0\in X$, set $k:=0$ and $\hat{\mathbb{N}}=\mathbb{N}\cup \{0\}$.
\item \textit{Step 2.} (Iterations)
 \begin{description}
\item $(i)$ If for $k\in \hat{\mathbb{N}}$, $y_k$ be a solution of $P_k$ : = $\left\{ \begin{array}{ll} \min \sum\limits_{i\in I}f_i(x), &\\
 \text{s. t. }X_k:=\{x\in X: f_i(x)\leq f_i(x_k),i\in I\}\end{array}\right.$ \\\textit{i.e.}, $y_k\in Sol(P_k)$ and  $y_k=x_k$, then stop. 
\item $(ii)$ Else, set $x_{k+1}:=y_k$ and $k:=k+1$. Go to step $2(i)$,
\end{description}
\end{description}
The idea of the algorithm can be implemented in the following way. Consider a tolerance factor $\epsilon_k\ge 0$ for $k\ge 1$ associated with the problem $P_k$. If $y_k$ solves $P_k$ and $y_k\ne x_k$ but $\|y_k-x_k\|\leq \epsilon_k$, then accept $x_k$ as the solution, otherwise set $x_{k+1}=y_k$ and solve $P_{k+1}$. As $k$ increases, we can keep on decreasing the tolerance limit $\epsilon_k$, \textit{i.e.}, $\epsilon_k\downarrow 0$ as $k\rightarrow \infty$. 

\begin{thm}\label{t31}
Let us consider the problem CVOP with each $f_i$ be a convex function and the feasible set $X$ be a compact set. Then,
\begin{itemize}
\item the iterates of the algorithm CC1 is well-defined.
\item If $\{x_k\}$ is a sequence generated by the above algorithm. Then any limit point of the sequence $\{x_k\}$ is a Pareto minimizer of the problem CVOP.
\end{itemize}
\end{thm}
\textit{Proof:} 
As $x_k\in X_k$, each $f_i$ is convex and the set $X_k$ is non empty, closed, convex and bounded, the set $Sol(P_k)$ is non-empty which proves the first part of the theorem.\\
If for some $k$, the solution $y_k$ of $P_k$ is equal to $x_k$, then one stops since by the Theorem \ref{c1tm34}, $x_k$ is a Pareto minimizer of the problem CVOP.\\
Now suppose for all $k\in \mathbb{N}$, the solution $y_k$ of $P_k$ is such that $y_k\ne x_k$. We shall now generate a sequence $\{x_k\}_{k\in \mathbb{N}}$, by defining $x_{k+1}=y_k$ for $k=0,1,\ldots$. Now since $X_k$'s are non-empty compact sets which are non-increasing,\textit{i.e.}, $X_{k+1}\subset X_k$, for all $k\in \mathbb{N}$, then by using Cantor's Intersection Theorem, we have $\bigcap\limits_{k= 1}^{\infty}X_k\not =\emptyset$. Further, as $x_k\in X_k$ for all $k\in \mathbb{N}$, $\{x_k\}$ is bounded and thus has a limit point, say $x^*$. Thus, there exists a subsequence $\{x_{k_j}\}_{j\in \mathbb{N}}$ such that $x_{k_j}\rightarrow x^*$. Further $x^*\in \bigcap\limits_{k= 1}^{\infty}X_{k_j} $, using the fact that $\bigcap\limits_{k= 1}^{\infty}X_{k_j} \ne \emptyset$\ as $X_{k_{j+1}}\subset X_{k_j}$ for all $j\in \mathbb{N}$. Hence, $x^*\in X_{k_j}$ for all $j\in \mathbb{N}$.
Now let us assume on the contrary that $x^*$ is not a Pareto minimizer, \textit{i.e.,} there exists $\hat{x}\in X$ such that 
\begin{equation}\label{eq31}
f_i(\hat{x})\leq f_i(x^*), \text{for all } i\in I \text{ and } f_r(\hat{x})<f_r(x^*), \text{for some } r\in I.
\end{equation} Hence,
\begin{equation}\label{eqN51}
\sum\limits_{i\in I}f_i(\hat{x})< \sum\limits_{i\in I}f_i({x}^*).
\end{equation}
Now as $x^*\in X_{k_j}$ for all $j\in \mathbb{N}$, $f_i(x^*)\leq f_i(x_{k_j})$ for all $j$. Thus,  \eqref{eq31} implies that $f_i(\hat{x})\leq f_i(x_{k_j})$ for all $j\ge 0$ and $i\in I$. Hence $\hat{x}\in X_k$ for all $k\in \mathbb{N}$. Further, $x_{k_j}\in Sol(P_{k_{j-1}})$, for $j\in \mathbb{N}$ which shows that
\begin{equation}\label{eqN52}
\sum\limits_{i\in I} f_i(x_{k_j})\leq\sum\limits_{i\in I} f_i(x)  \text{ for all } x\in X_{k_{j-1}}.
\end{equation}
As we have argued that $\hat{x}\in X_k$ for all $k\in \mathbb{N}$, from \eqref{eqN52} we have
\begin{equation}\label{eqN53}
\sum\limits_{i\in I} f_i(x_{k_j})\leq\sum\limits_{i\in I} f_i(\hat{x}) . 
\end{equation}
As $j\rightarrow \infty $ in \eqref{eqN53} and using the continuity of the functions $f_i$ for all $i\in I$, we get 
$$\sum\limits_{i\in I}f_i(x^*)\leq\sum\limits_{i\in I}f_i(\hat{x}),$$
which contradicts the inequality \eqref{eqN51}. Hence $x^*$ is a Pareto minimizer of the problem CVOP which completes the proof of the second part of the theorem.\hfill$\Box$
\begin{remark}\rm
We can generalize the Theorem \ref{t31} by replacing the assumption of convexity of each function $f_i$ by proper, lower semicontinuous and $\mathbb{R}^m_+$-convex objective function $f$ and compactness of the feasible set $X$ by $\mathbb{R}^m_+$-completeness of the set $(f(x_0)-\mathbb{R}^m_+)\cap f(X) $ which means that for all sequences $\{a_n\}\subset X$ with $a_0=x_0$ such that $f(a_{n+1})\leq f(a_n)$ for all $n\in \mathbb{N}$, there exists $a\in X$ such that $f(a)\leq f(a_n)$ for all $n\in \mathbb{N}$ (see \cite{bonnel2005proximal}).
\end{remark}
Let us again take a careful look at the algorithm CC1. We don not build any separate subroutine for solving $P_k$ but use the CVX tool box in MATLAB. Given the problem $P_k$ , one of the key ideas here is to see how near the solution of $P_k$ is to $x_k$. What may happen that one might have $x_k\in Sol(P_k)$ and CVX provides a solution $y_k\ne x_k$. Then we have missed the Pareto minimizer $x^k$ and the run time of the algorithm gets increased. How does one address this particular issue?\
An intuitive idea that comes into mind is that for each $k\in\hat{\mathbb{N}}$, one can try to see what will happen if we choose $x_k$ as the starting solution for $P_k$? Will the algorithm CVX return $x_k$ as the solution if $x_k$ is indeed a member of $Sol(P_k)$ or at least it gives us an $y_k$ with $\|y_k-x_k\|\le \epsilon_k$, where $\epsilon_k>0$ is a threshold value for rejecting $x_k$ as the solution of $P_k$. We did some numerical experiments and found that at least in those cases when $x_k\in Sol(P_k) $ and we consider $x_k$ as our starting solution, the algorithm CVX indeed returned $x_k$ as the solution.

Another route often suggested under such circumstances is to convert $P_k$ into a strongly convex program so that CVX will give us just a unique solution for each $k\in \hat{\mathbb{N}}$. In fact let us consider the following scalar problem $\hat{P}(f,X,x^*)$ given as 
\begin{equation*} 
 \left\{ \begin{array}{ll} \min \sum\limits_{i\in I}(f_i(x)+\frac{1}{2}\|x-x^*\|^2), & \\
{ \text{subject to }}~~ f_i(x)-f_i(x_0)\leq 0, &\text{ for all }  i\in I, \\  
x\in X. & \end{array}\right.
\end{equation*}
It can be easily proved that $x^*$ is a Pareto minimizer of CVOP if and only if $x^*$ is a solution of the scalar strongly convex program $\hat{P}(f,X,x^*)$. In fact in the algorithm CC1, we can replace $P_k$ with $\hat{P_k}$ given as
\begin{equation*} 
 \left\{ \begin{array}{ll} \min \sum\limits_{i\in I}(f_i(x)+\frac{1}{2}\|x-x^*\|^2), & \\ 
x\in X_k. & \end{array}\right.
\end{equation*}
However, having a strongly convex problem $P_k$ does not completely alleviate the issue at hand. If we are looking for a solution with very high accuracy, then for each $k$, our threshold limit is very small and even for the problem $\hat{P_k}$, the routine CVX might return $y_k$ as the solution where $\|y_k-x_k\|> \epsilon_k$. Thus, we believe that we can just work with $P_k$ and the numerical experiments given below showed us that the algorithm CC1 is working efficiently. The use of scalar problems similar to $\hat{P}(f,X,x^*)$ does not appear to be new. In the study of proximal point method for vector optimization problem, Bonnel et al. \cite{bonnel2005proximal} use the following scalar problem at $k$-th step of their algorithm denoted by ALG1:
\begin{equation*} 
 \left\{ \begin{array}{ll} \min ~\langle f(x)+\frac{\alpha_k}{2}\|x-x_k\|^2e_k,z_k\rangle, & \\
{ \text{subject to }}~~ f_i(x)-f_i(x_k)\leq 0, &\text{ for all }  i\in I, \\  
x\in X, & \end{array}\right.
\end{equation*}
where $\{z_k\}$, $\{e_k\}$ and $\{\alpha_k\}$ are some sequences in the considered ordering cone $C\subset \mathbb{R}^m$ and in our case $C=\mathbb{R}^m_+$. There are two major differences between our algorithms CC1 and the algorithm AGL1 of \cite{bonnel2005proximal}. Firstly, ALG1 needs a new set of parameters such as $z_k$, $e_k$ and $\alpha_k$ at every step of the algorithm whereas CC1 does not need any such parameters. Secondly, the stopping criteria of ALG1 involve solving the multiobjective problem which is not a practical approach when it comes to implementing the algorithm. On the other hand, stopping criteria of CC1 is fairly simple which only involve checking the error is within the considered threshold limit. Further, Bonnel et al. \cite{bonnel2005proximal} admitted that it is not easy to implement their algorithm to numerical test problems whereas the algorithm CC1 is easily applicable to the test problems using MATLAB software. Similar to Bonnel et al. \cite{bonnel2005proximal}, there are several point-by-point algorithms to solve vector optimization problems. Some algorithms can only guarantee that the resultant solution is a weak Pareto minimizer or critical point, for example, \cite{drummond2004projected}, \cite{fukuda2011convergence}, whereas some algorithms need lots of assumption on objective functions which narrows down the applicability of the algorithm for applications, for example, \cite{fliege2009newton}, \cite{fliege2000steepest}. 
\section{Numerical Computations}
 In this section, we consider some convex vector optimization test problems to show the efficiency of the proposed algorithm CC1. The design of the algorithm CC1 has a advantage which makes it more useful than evolutionary algorithms, for example, NSGA-II \cite{deb2002fast} used to solve CVOP. Firstly, algorithm CC1 has a convergence analysis which assures that the solution will be a Pareto solution. On the other hand, NSGA-II does not have any convergence analysis. Secondly, we can trace Pareto solution corresponding to the efficient frontier of the problem in the algorithm CC1 while NSGA-II only gives the information of the efficient frontier. But we need to admit that the run time of algorithm CC1 is higher than NSGA-II.\\~
 
\noindent In each of the examples considered below, we 
 consider convex bi-objective test problems where the objective functions $f_1$ and $f_2$ and the feasible set $X$ are separately defined. All problems are solved for 200 randomly generated starting points between the lower bound and upper bound of the decision variables of the considered problem. At each step, we use CVX software on the MATLAB platform to solve the scalar problem $P_k$ of the algorithm CC1. In figures,``+" denotes the image of starting point and ``o" denotes the image of the resultant point of the algorithm CC1 under the objective function $f=(f_1,f_2)$. We tested algorithm CC1 on convex vector optimization problems with smooth as well as non-smooth input data to check the efficiency of this algorithm. First three examples deals with convex vector optimization with smooth objective functions and the objective functions of the last example is non-smooth. For simplicity, we consider only two dimensional problems but the algorithm CC1 is applicable for higher dimensional problems as well.\\
 
 \begin{example} \label{ex1}
 Let $f_1,f_2:\mathbb{R}\to \mathbb{R}$ be defined by $f_1(x)=x^2,~f_2(x)=(x-2)^2$ and let $X=\{x\in \mathbb{R}:-10\leq x\leq 10\}$ (Schaffer problem \cite{schaffer1986some}).
 \begin{figure}[H]
\begin{center}
\includegraphics[width=0.62\textwidth]{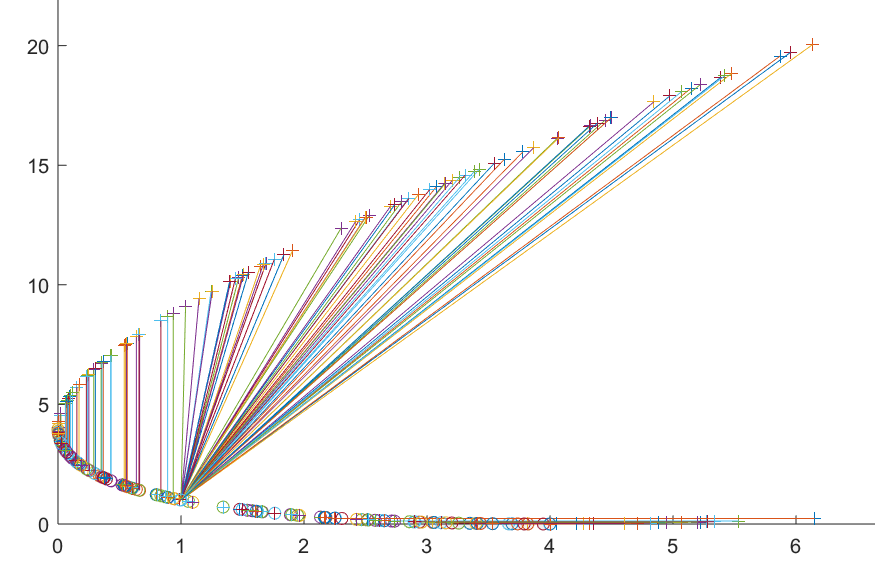}\caption{ Efficient frontier of Example \ref{ex1}}
\label{fig:4.5}
\end{center}
\end{figure}
 \end{example}
 \begin{example}\label{ex2}
  Let $f_1,f_2:\mathbb{R}^2\to \mathbb{R}$ be defined by $f_1(x_1,x_2)=x_1,~f_2(x_1,x_2)=x_2$ and let $X=\{(x_1,x_2)\in \mathbb{R}^2:x_1^2+x_2^2\leq 1\}$ (Jahn problem \cite{JAHN}). 
  \begin{figure}[H]
\begin{center}
\includegraphics[width=0.62\textwidth]{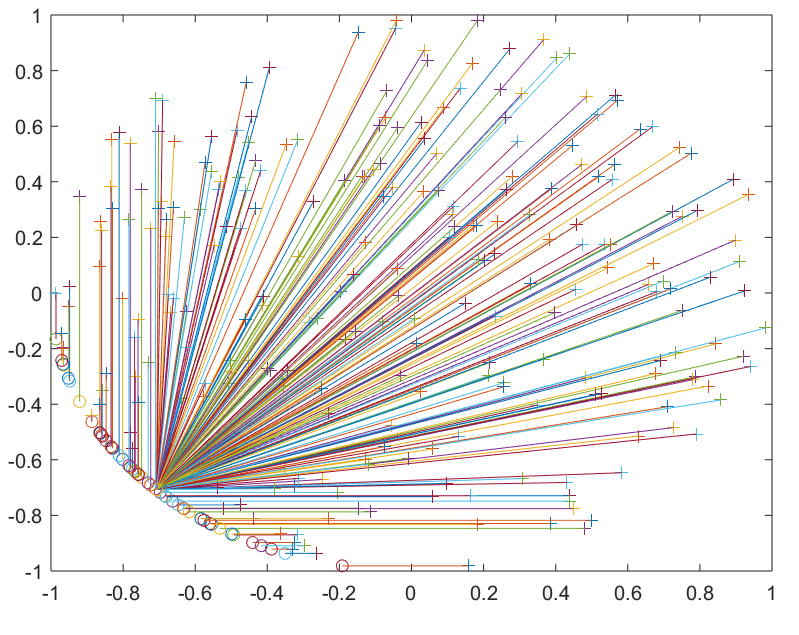}\caption{ Efficient frontier of Example \ref{ex2} }
\label{fig:4.5}
\end{center}
\end{figure}
 \end{example}
\begin{example}\label{ex3}
  Let $f_1,f_2:\mathbb{R}^2\to \mathbb{R}$ be defined by $f_1(x_1,x_2)=4x_1^2+4x_2^2,~
f_2(x_1,x_2)=(x_1-5)^2+(x_2-5)^2$ and $X=\{(x_1,x_2)\in \mathbb{R}^2: 
-5\le x_1\le 10, -5\le x_2\le 10\}$ (Binh problem \cite{binh1999multiobjective}). 
  \begin{figure}[H]
\begin{center}
\includegraphics[width=0.59\textwidth]{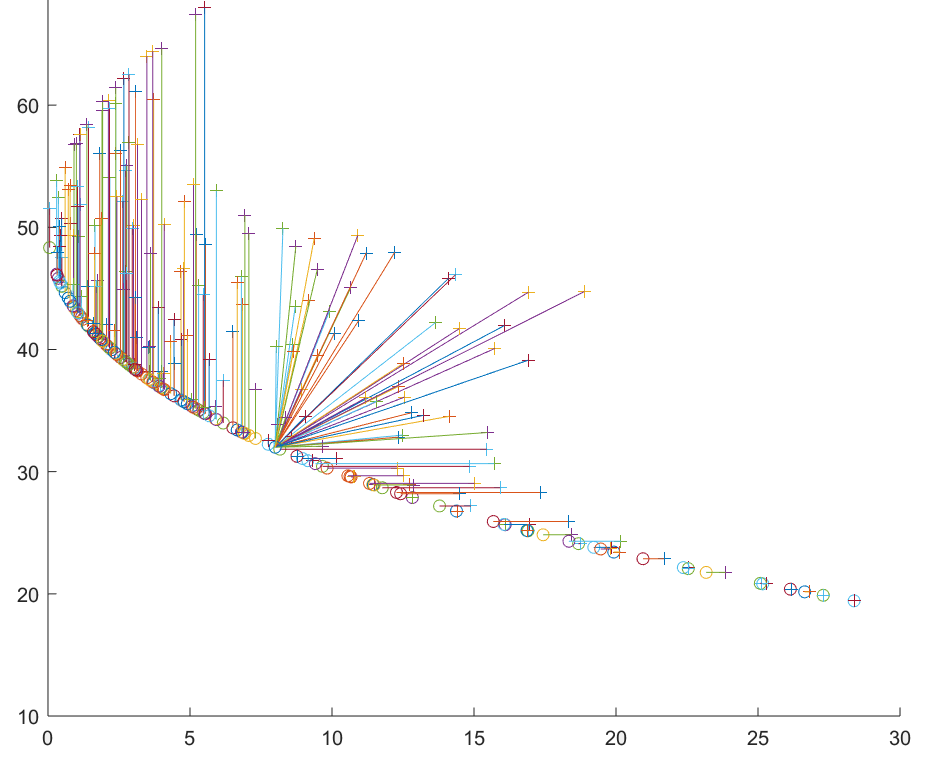}\caption{ Efficient frontier of Example \ref{ex3}}
\label{fig:4.5}
\end{center}
\end{figure}
 \end{example}
 
\begin{example}\label{ex4}
  Let $f_1,f_2:\mathbb{R}^2\to \mathbb{R}$ be defined by $f_1(x_1,x_2)=\max \{x_1,x_2\},~
f_2(x_1,x_2)=|x_1|+|x_2|$ and $X=\{(x_1,x_2)\in \mathbb{R}^2: 
-2\le x_1\le 2, -2\le x_2\le 2\}$. 
  \begin{figure}[H]
\begin{center}
\includegraphics[width=0.62\textwidth]{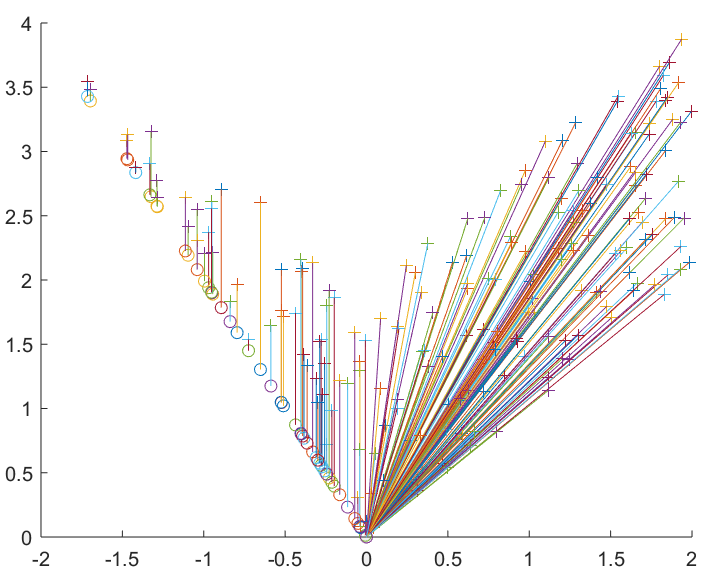}\caption{ Efficient frontier of Example \ref{ex4}}
\label{fig:4.5}
\end{center}
\end{figure}
 \end{example}

\end{document}